\documentclass[11pt]{article}
\usepackage{amsfonts,latexsym,amsmath,amscd,geometry}
\usepackage{color}
\usepackage{amssymb}
\usepackage{latexsym}
\usepackage{xcolor}
\usepackage[title]{appendix}

\newcommand \nc{\newcommand}
\newtheorem{theorem}{Theorem}[section]
\newtheorem{lemma}[theorem]{Lemma}

\nc{\ba}{\begin{array}}\nc{\ea}{\end{array}}
\nc{\be}{\begin{eqnarray}}\nc{\ee}{\end{eqnarray}}
\nc{\beq}{\begin{equation}}\nc{\eeq}{\end{equation}}
\nc{\bex}{\begin{eqnarray*}}\nc{\eex}{\end{eqnarray*}}
\nc{\btm}{\begin{theorem}} \nc{\etm}{\end{theorem}}
\nc{\blm}{\begin{lemma}} \nc{\elm}{\end{lemma}}
\nc{\R}{\mathbb{R}}  
\nc{\e}{\mathbf{e}} 
\nc{\va}{\varphi}
\nc{\ve}{\varepsilon}

\def\pf{\noindent{\bf Proof.\quad}}
\def\endpf{\hfill$\Box$}
\def\di{\mbox{div\,}}

\begin{document}

\title{On singularities of Ericksen-Leslie system in dimension three}

\author{Tao Huang\footnote{Department of Mathematics, Wayne State University, Detroit, MI 48202, USA.}  \quad Peiyong Wang$^{*}$\\
}
\date{}
\maketitle

\begin{abstract}
In this paper, we consider the initial and boundary value problem of Ericksen-Leslie system modeling nematic liquid crystal flows in dimension three. Two examples of singularity at finite time are constructed. The first example is constructed in a special axisymmetric class with suitable axisymmetric initial and boundary data, while the second example is constructed for an initial data with small energy but nontrivial topology. A counter example of maximum principle to the system is constructed by utilizing the Poiseuille flow in dimension one.

\end{abstract}

\section{Introduction}

Nematic liquid crystals are composed of rod-like molecules characterized by the average alignment of  the long axes of neighboring molecules, which have simplest structure among various types of liquid crystals. The dynamic theory of nematic liquid crystals has been first proposed by Ericksen \cite{ericksen62} and Leslie \cite{leslie68} in the 1960's, which is a macroscopic continuum description of the time evolution of both the flow velocity field and the orientation order parameter of rod-like liquid crystals.
More precisely, we would consider the following Ericksen-Leslie system in $\Omega\times(0,\infty)$, where $\Omega\subset\R^3$ is a bounded domain with smooth boundary
\begin{equation}\label{EL}
\begin{cases}
u_t+u\cdot\nabla u+\nabla P=-\nabla\cdot\big(\nabla d\odot\nabla d\big)+\nabla \cdot(\sigma^L(u,d)),\\
\nabla \cdot u=0, \\
\lambda_1(d_t+u\cdot\nabla d-\Lambda d)+\lambda_2Ad=\Delta d+|\nabla d|^2d+\lambda_2(d^TAd)d,
\end{cases}
\end{equation}
where 
$u(x,t): \Omega\times(0,\infty)\rightarrow \R^3$ is the fluid velocity field, $d(x,t):\Omega\times(0,\infty)\rightarrow \mathbb S^2$ is the orientation order parameter of nematic material at $(x,t)$, and $P(x,t): \Omega\times(0,\infty)\rightarrow \R$ is the pressure. Denote $u=(u_1, u_2, u_3)\in  \R^3$ and $d=(d_1, d_2, d_3)\in\mathbb S^2$. Then
$$
\big(\nabla d\odot\nabla d\big)_{ij}
=\nabla_i d\cdot \nabla_j d,
$$
and 
$$
A_{ij}=\frac{1}{2}\left(\frac{\partial u_j}{\partial x_i}+\frac{\partial u_i}{\partial x_j}\right),\quad \Lambda_{ij}=\frac{1}{2}\left(\frac{\partial u_i}{\partial x_j}-\frac{\partial u_j}{\partial x_i}\right),\quad 
N_i=\partial_td_i+u\cdot\nabla d_i-\Lambda_{ij}d_j.
$$
denote the rate of the strain tensor, the skew-symmetric part of the strain rate, and the
rigid rotation part of the direction changing rate by fluid vorticity, respectively.  The left side of the third equation in \eqref{EL} is the kinematic transport, which represents the effect of the macroscopic flow field on the microscopic structure. The material coefficients $\lambda_1$ and $\lambda_2$ reflect the molecular shape and the slippery part between the fluid and the particles. The term with $\lambda_1$ represents the rigid rotation of molecules, while the term with $\lambda_2$ stands for the stretching of molecules by the flow. The viscous (Leslie) stress tensor $\sigma^L$ has the following form (cf. \cite{Les79})
$$
\sigma_{ij}^L(u,d)=\mu_1 d_kd_pA_{kp}d_id_j+\mu_2 N_id_j+\mu_3d_iN_j+\mu_4A_{ij}+\mu_5A_{ik}d_kd_j+\mu_6 A_{jk}d_kd_i.
$$
The viscous coefficients $\mu_i$, $i=1,\cdots, 6,$ are called the Leslie's coefficients and the following relations are often assumed in the literature
\beq\label{parodi}
\mu_2+\mu_3=\mu_6-\mu_5
\eeq
\beq\label{nec1}
\lambda_1=\mu_3-\mu_2>0, \quad \lambda_2=\mu_6-\mu_5
\eeq
%\beq\label{nec2}
%\mu_1+\frac{\lambda_2^2}{\lambda_1}\geq 0,\quad 
%\mu_4>0, \quad \mu_5+\mu_6\geq \frac{\lambda_2^2}{\lambda_1}.
%\eeq
\beq\label{mus}
\mu_4>0,\quad 2\mu_1+3\mu_4+2\mu_5+2\mu_6>0,
\quad 2\mu_4+\mu_5+\mu_6>\frac{\lambda_2^2}{\lambda_1}.
\eeq
The first relation is called the  Parodi's relation, which has been derived from the Onsager reciprocal relations expressing the equality of certain relations between flows and forces in thermodynamic systems out of equilibrium (cf. \cite{Parodi70}).
The second set of relations are the compatibility conditions. The third empirical relations are necessary to obtain the energy inequality (cf. \cite{Les79}, \cite{WZZ13}). Throughout the paper, we assume that \eqref{parodi}, \eqref{nec1} and \eqref{mus} are valid for system \eqref{EL}.

We consider the initial data 
\beq\label{initial}
u(x,0)=u_0,\quad d(x,0)=d_0
\eeq
with 
%$u_0\in L^2(\Omega,\R^3)$, $d_0\in H^1(\Omega, \mathbb S^2)$ 
$\di u_0=0$ and $|d_0|=1$, 
and the boundary data 
\beq\label{bdry}
u(x,t)|_{\partial \Omega\times[0,\infty)}=0,\quad d(x,t)|_{\partial \Omega\times[0,\infty)}=\e
\eeq
with compatible conditions
\beq
u_0|_{\partial\Omega}=0,\quad d_0|_{\partial\Omega}=\e.
\eeq
Here $\e=(0,0,1)$ is a constant vector.

In dimension two, the existence of global weak solutions to the Cauchy problem of \eqref{EL} has been established in \cite{huanglinwang14}. The weak solution has been proved to be regular except for finitely many times (see also \cite{hongxin12} for related results). The uniqueness of such a weak solution has been proved in \cite{LTX16, wangwang14}. We would like to point out that the assumptions of Leslie coefficients in \cite{huanglinwang14} is stronger than ours, \eqref{mus}. However, their results are still valid with weaker assumption \eqref{mus} since we can prove the similar energy inequality (see Lemma \ref{lemma0} below for more details). 
In dimension three, the global well-posedness combining with long time behaviors for the system \eqref{EL} around equilibrium under various assumptions on the Leslie coefficients has been studied in \cite{WZZ13, WXL13, hieberpruss17}.

There is a simplified system that has been first proposed in \cite{lin89} by neglecting the Leslie stress. There have been many results on the existence and partial regularity of this simplified system in \cite{linlinwang10, linwang10, linwang16, HLLW16, llwwz19} and the references therein.

\medskip
The paper is organized as follows. In Section 2, we will state our main results. In Section 3, a special form of axisymmetric solution will be derived. In Section 4, global existence and singularities at finite time will be discussed.  In Section 5, another example of finite time singularity will be constructed with an initial data with small energy but with large topology. In Section 6, we will construct a counter example to show that the maximum principle is not necessarily valid for Ericksen-Leslie system via the Poiseuille flow in dimension one.

\section{Main results}
 \setcounter{equation}{0}
\setcounter{theorem}{0}

Let $B_1^n\subset\mathbb R^n$ denote the unit ball centered at $0$. Inspired by the results in \cite{HLLW16}, we first consider the domain $\Omega=B_1^2\times [0,1]$ 
and an axisymmetric solution $(u,P, d)$ to the Ericksen-Leslie system \eqref{EL} in the following special form
\beq\label{specialsol}
\begin{cases}
u(r,\theta, z,t):=v(r,t)\e^r+w(z,t)\e^3,\\
d(r,\theta, z, t):=\sin\va(r,t) \e^r+\cos\va(r,t) \e^3,\\
P(r,\theta, z,t):=Q(r,t)+R(z,t).
\end{cases}
\eeq
For any ${\bf x}=(x,y,z)\in\Omega$, we take the initial data as follows
\beq\label{initial1}
u_0({\bf x})=(x,y,-2z),
\eeq
and
\beq\label{initial2}
\displaystyle d_0({\bf x})=\Big(\frac{x}{\sqrt{x^2+y^2}}\sin\varphi_0\big(\sqrt{x^2+y^2}\big),
\frac{y}{\sqrt{x^2+y^2}}\sin\varphi_0\big(\sqrt{x^2+y^2}\big), \cos\varphi_0\big(\sqrt{x^2+y^2}\big)\Big),
\eeq
with $\varphi_0\in C^\infty([0,1])$ and 
$\varphi_0(0)=0$. 
To derive the special axisymmetric form of system \eqref{EL}, we also need to consider the following boundary conditions
%\begin{equation}
%\big(u({\bf x}, 0),d({\bf x}, 0)\big)=\big(u_0({\bf x}), d_0({\bf x})\big),\ {\bf x}\in\Omega, \label{initial1}
%\end{equation}
\begin{equation}\label{boundary}
\begin{cases}
\ u({\bf x},t)=u_0({\bf x}) &\ {\bf x}\in\partial\Omega, \ t>0,\\
\ d({\bf x}, t)=d_0({\bf x}) & \ {\bf x}\in \partial B_1^2\times [0,1], \ t>0,\\
\frac{\partial d}{\partial z} ({\bf x}, t)=0 & \ {\bf x}\in B_1^2\times \{0, 1\}, \ t>0.
\end{cases}
\end{equation}
Hence the Ericksen-Leslie system \eqref{EL} becomes (see Section \ref{secderaxi} for details)
$$v(r,t)=r,\quad w(z,t)=-2z,$$
\beq\label{ELaxeqn0}
\lambda_1(\va_t+r\va_r)=\va_{rr}+\frac{1}{r}\va_r-\frac{\sin(2\va)}{2r^2}-3\lambda_2\sin\va\cos\va
\eeq

Our first result concerns the global existence of the equation \eqref{ELaxeqn0}, which implies the global existence of the Ericksen-Leslie system \eqref{EL}.
\btm\label{glbth1}{\it Suppose $\phi_0\in C^\infty([0,1])$
satisfies $\va_0(0)=0$ and
$$0\leq \va_0(r)< \pi,$$ 
for all $r\in [0,1]$.
Then there is no finite time blowup for the smooth solution 
to \eqref{ELaxeqn0}.
}
\etm

We can also construction an example of singularity to the system \eqref{EL} in this axisymmetric class.

\btm\label{blowthm2}{\it There exists $\phi_0\in C^\infty([0,1])$, with $\phi_0(0)=0$
and $\va_0(1)> \pi$, such that
the short time smooth solution $\va$ to \eqref{ELaxeqn0},
with initial and boundary conditions \eqref{initial1}-\eqref{boundary},
must blow up at $T_0$ for some
$0< T_0=T_0(\phi_0)<+\infty$. More precisely,  $\va_r(0,t)\rightarrow\infty$ as $t\rightarrow T_0^-$.
}
\etm

The main idea to prove the Theorems \ref{glbth1} and \ref{blowthm2} is to construct suitable supersolution and subsolution, which has been first utilized in the breakthrough work \cite{chang-ding-ye} when the authors studied the singularity of harmonic heat flows. Similar results have been proved for simplified Ericksen-Leslie system in \cite{HLLW16}. However, there are still two difficulties for the full Ericksen-Leslie system \eqref{EL}. The first one is that the comparison principle and maximum principle may not be valid for the system \eqref{EL} because of the terms related to the skew-symmetric part $\Lambda$ of the strain rate and the symmetric part $A$. Fortunately, it holds $\Lambda=0$, and $A$ is a diagonal matrix for our special axisymmetric solution \eqref{specialsol} (cf. Section \ref{secderaxi}). Therefore, we still have comparison and maximum principles for equation \eqref{ELaxeqn0} (cf. Section \ref{secapri}). The second difficulty comes from the nonlinear term $-3\lambda_2\sin\va\cos\va$ in \eqref{ELaxeqn0}. It is not hard to see that this term will not produce any singularity. However, since we don't know the sign of $\lambda_2$, it brings more technique difficulties to the construction of a supersolution and a subsolution.

We also want to point out that the special class of axisymmetric solutions may not satisfy the following energy inequality, which has been proved in \cite{WZZ13} and \cite{AnnaLiu19}. 
\begin{lemma}\label{lemma0}
Suppose that $(u,d)$ is a regular solution to the system \eqref{EL} with the initial and boundary conditions \eqref{initial} and \eqref{bdry} under the assumptions \eqref{parodi}, \eqref{nec1} and \eqref{mus}. For any $t\in(0,\infty)$, the following energy equality holds
\beq\label{enginq}
\frac{d}{dt}\frac{1}{2}\int_{\Omega}\left(|u|^2+|\nabla d|^2\right)\,
\leq-\int_{\Omega}\left(\alpha_0|\nabla u|^2+\frac{1}{\lambda_1}\left|\Delta d+|\nabla d|^2d\right|^2\right)\,,
\eeq
for some constant $\alpha_0>0$.
\end{lemma}
Therefore, we would like to construct another example of finite time singularity of (\ref{EL}) for a generic
initial data, in which the solution satisfies the energy dissipation inequality.  
For simplicity, we take $\Omega=B_1^3$, the unit ball in dimension three, and denote 
$$C^\infty_{0,{\rm{div}}} (B^3_1,\mathbb R^3) 
:= \Big\{v \in C^\infty(B_1^3,\mathbb R^3) \ \big| \ \nabla\cdot v =  0\Big\},$$
$$C_{\bf e}^\infty(B_1^3,\mathbb S^2):=\Big\{d\in C^\infty(B_1^3,\mathbb S^2)\ \big| \
d={\bf e}\ {\rm{on}}\ \partial B_1^3\Big\}.$$
For continuous maps $f, g\in C\big(\overline{B_1^3},\mathbb S^2\big)$, with $f=g$ on $\partial B_1^3$,
we say that $f$ is homotopic to $g$ relative to $\partial B_1^3$ if there exists a continuous map
$\Phi\in C\big(\overline{B_1^3}\times [0,1],\mathbb S^2\big)$ such that $\Phi(\cdot, t)=f(\cdot)=g(\cdot)$ on $\partial B_1^3$, for all $0\le t\le 1$, and
$\Phi(\cdot, 0)=f(\cdot)$ and $\Phi(\cdot, 1)=g(\cdot)$ in $B_1^3$. Now we are ready to state our main result.

\begin{theorem} \label{blowthm1} There exists 
$\epsilon_0 >0$ such that if $u_0 \in C^\infty_{0, {\rm{div}}}(B_1^3, \mathbb R^3)$ and 
$d_0 \in C_{\bf e}^\infty(B_1^3, \mathbb S^2)$ satisfy that
$d_0$ is not homotopic to the constant map ${\bf e}:B_1^3\to\mathbb S^2$ relative to $\partial B_1^3$, and
\begin{equation}\label{seng}
E(u_0,d_0):=\frac12\int_{B_1^3}(|u_0|^2+|\nabla d_0|^2)\le\epsilon_0^2.
\end{equation}
Then the short time smooth solution $(u, d, P) : B_1^3 \times [0, T) \to \mathbb R^3 \times\mathbb S^2\times \R$
to the Ericksen-Leslie system (\ref{EL}),
with the initial-boundary condition \eqref{initial} and \eqref{bdry} under the assumptions \eqref{parodi}, \eqref{nec1} and \eqref{mus}
%\begin{equation}\label{IBC}
%\begin{cases}
%(u,d)\big|_{t=0} = (u_0,d_0),\ {\rm{in}}\  B_1^3,\\
%(u,d)\big|_{\partial B_1^3} = (0,{\bf e}), \ 0 < t < T,
%\end{cases}
%\end{equation}
 must blow up before time $T = 1$.
\end{theorem}

Due to the energy estimate \eqref{enginq}, the proof of Theorem \ref{blowthm1} is quite similar to those in \cite{HLLW16} for the simplified system. For completeness, we will briefly sketch the main steps of the proof in Sec \ref{sec2ex}.

We would like to remark that there does exist $(u_0,d_0)\in C^\infty_{0, {\rm{div}}}(B_1^3,\mathbb R^3)
\times C^\infty_{\bf e}(B_1^3, \mathbb S^2)$ satisfying the conditions of Theorem \ref{blowthm1} which has been previously used by \cite{HLLW16} in construction of finite time singularity of simplified liquid crystal flows (see also \cite{linwangbook}, \cite{DW} for other applications). 
More precisely, let $H(z,w)=(|z|^2-|w|^2, 2zw):\mathbb S^3\equiv \big\{(z,w)\in \mathbb C\times\mathbb C: |z|^2+|w|^2=1\big\}\to \mathbb S^2\subset \mathbb R\times \mathbb C$ be the Hopf map. Let $D_\lambda({\bf x})=\lambda {\bf x}:\mathbb R^3\to\mathbb R^3$ be the dilation map for $\lambda>0$, $\Pi:\mathbb S^3\to\overline{\mathbb R^3}$ be the 
stereographic projection map from ${\bf e}$, and $\Psi_\lambda=\Pi^{-1}\circ D_\lambda\circ \Pi:\mathbb S^3\to\mathbb S^3$. Then direct calculations imply that the Dirichlet energy of $H\circ\Psi_\lambda:\mathbb S^3\to\mathbb S^2$
satisfies
$$\lim_{\lambda\rightarrow\infty}\int_{\mathbb S^3}\big|\nabla(H\circ \Psi_\lambda)\big|^2\,d\sigma=0.$$
Moreover, it is easy to see that $H\circ\Psi_\lambda$ is not homotopic to the constant map ${\bf e}:\mathbb S^3\to\mathbb S^2$. Let $\Phi\in C^\infty\big(\overline{B_1^3}, \mathbb S^3\big)$ such that $\Phi: B_1^3\to \mathbb S^3\setminus\{\bf e\}$ is a diffeomorphism and $\Phi={\bf e}$ on $\partial B_1^3$. Now we can check that for any
$u\in C^\infty_{0,{\rm{div}}}(B_1^3,\mathbb R^3)$, since
$$\lim_{\lambda\rightarrow \infty} E(\lambda^{-1} u, H\circ\Psi_\lambda\circ\Phi)=0,$$
we can find a sufficiently large $\lambda_0>0$ depending on $u$, $H$, and $\Phi$ such that
$$(u_0,d_0):=(\lambda_0^{-1} u, H\circ\Psi_{\lambda_0}\circ\Phi):B_1^3\to\mathbb R^3\times\mathbb S^2$$
satisfies the condition (\ref{seng}) of Theorem \ref{blowthm1}, and $d_0$ is not homotopic to the constant map
${\bf e}$ relative to $\partial B_1^3$.

It is also an interesting topic to investigate more on the homotopic condition of $d_0$ in Theorem \ref{blowthm1}. If we assume that 
$d_0$ is homotopic to the constant map ${\bf e}$ relative to $\partial B_1^3$, then the angle $\va$ between $d_0$ and ${\bf e}$ should belong to $[0,\pi)$. Because the maximum principle (cf. Lemma \ref{maxp}) may not be valid for the system \eqref{EL}, it is natural to ask if $d_0$ is homotopic to the constant map ${\bf e}$ relative to $\partial B_1^3$ and the initial energy of the system is small, whether the angle $\va$ between $d_0$ and ${\bf e}$ can go beyond $\pi$ at some time $t$. If this happens, $d(\cdot, t)$ is not homotopic to the constant map ${\bf e}$ relative to $\partial B_1^3$, which, combined with the energy inequality \eqref{enginq} and Theorem \ref{blowthm1}, implies the finite time singularity. Unfortunately, we can't verify this argument in dimension three due to technical difficulties. However, utilizing the Poiseuille flows in dimension one, we can construct a counter example of the maximum principle (cf. Lemma \ref{maxp}), which does imply the maximum principle may not be valid in the general case of Ericksen-Leslie system \eqref{EL} (see Section \ref{seccnex} for more details).

\section{Axisymmetric solutions}
 \setcounter{equation}{0}
\setcounter{theorem}{0}

Let $(r,\theta, z)$ denote the cylindrical
coordinates of $\mathbb R^3$, and set
$$\e^r=(\cos\theta,\sin\theta,0)^T,\quad \e^{\theta}=(-\sin\theta,\cos\theta,0)^T,\quad \e^3=(0,0,1)^T$$
as the canonical orthonormal base of $\mathbb R^3$ in the cylindrical coordinates.
For $\alpha\in [0,2\pi]$, let $R_\alpha\in {\rm{SO}}(3)$ denote the rotation map 
of angle $\alpha$ with respect to the $z$-axis. 
Recall that a vector field $v:\mathbb R^3\to\mathbb R^3$ is axisymmetric if
$$R_\alpha^{-1}\circ v\circ R_\alpha =v, \ \forall\ \alpha\in [0,2\pi].$$
Hence any axisymmetric vector field $v$ can be written as
\beq\label{axis1.1}
v(r,\theta,z)=v^r(r, z)\e^r+v^{\theta}(r,z)\e^{\theta}+v^3(r,z)\e^3.
\eeq
If, in addition, $v^\theta\equiv 0$, we say $v$ is axisymmetric without swirls.
A solution $(u,d, P)$ to the Ericksen-Leslie system (\ref{EL}) is said to be axisymmetric without swirls, if it holds that
$$\begin{cases} u(r,\theta, z,t)=u^{r}(r,z,t)\e^r+u^{3}(r,z,t)\e^3,\\
d(r,\theta, z,t)=\sin\va(r,z,t) \e^r+\cos\va(r,z,t) \e^3,\\
P(r,\theta,z,t)=P(r,z,t).
\end{cases}
$$
A domain $\Omega\subset\mathbb R^3$ is axisymmetric if it is invariant under a rotation
map $R_\alpha$ for any $\alpha\in [0,2\pi]$.

\subsection{Derivation of a special axisymmetric system}
\label{secderaxi}

It is not hard to see that the initial condition \eqref{initial1} for $(v,w,\phi)$ reduces to
\begin{equation}\label{sbinitial}
\begin{cases}
v|_{t=0}=r,\ \ \ \ \ \ \ 0\le r\le 1,\\
w|_{t=0}=-2z, \ \ 0\le z\le 1,\\
\va|_{t=0}=\va_0(r),\ 0\le r\le 1,
\end{cases}
\end{equation}
for some $\va_0\in C^\infty([0,1])$, with $\va_0(0)=0$.
The boundary condition \eqref{boundary} becomes
\beq\label{sboundary}
\begin{split}
\begin{cases} v(0,t)=0\\
v(1,t)=1,
\end{cases} & \quad \begin{cases} w(0,t)=0\\ w(1,t)=-2,\end{cases}
\quad \begin{cases} \va(0,t)=0\\ \va(1,t)=\va_0(1).\end{cases}
\end{split}
\eeq

\begin{lemma} \label{static} For $0<T\le +\infty$, suppose that $v, w\in C^\infty([0,1]\times [0, T))$ satisfies 
\begin{equation}\label{div}
\displaystyle\frac{1}{r}(rv)_r+w_z=0, \ (r,z)\in [0,1]\times [0,1],
\end{equation}
and
\beq\label{2pt}
\begin{split}
\begin{cases} v(0,t)=0\\
v(1,t)=1,
\end{cases} & \quad \begin{cases} w(0,t)=0\\ w(1,t)=-2.\end{cases}
\end{split}
\eeq
Then $v(r,t)=r$ for any $(r,t)\in [0,1]\times [0,T)$, and $w(z,t)=-2z$ for any $(z,t)\in [0,1]\times [0,T)$.
\end{lemma} 
\pf  Differentiating \eqref{div} with respect to $z$ yields
\bex
w_{zz}(z,t)=0
\eex
so that $w(z,t)=a_1(t)z+a_2(t)$ for some functions $a_1(t)$ and $a_2(t)$. 
Since $w(0,t)=0$ and $w(1,t)=-2$, we see that
$a_2(t)\equiv 0$ and $a_1(t)\equiv -2$. Thus $w(z,t)=-2z.$

Similarly, differentiating \eqref{div} with respect to $r$ yields
\bex
\big(\frac{1}{r}(rv)_r\big)_r(r,t)=0,
\eex
which implies that $rv(r,t)=b_1(t)r^2+b_2(t)$ for some functions $b_1(t)$ and $b_2(t)$. 
Since $v(0,t)=0$ and $v(1,t)=1$, we  see that
$b_2(t)\equiv 0$ and $b_1(t)\equiv 1$. Thus $v(r,t)=r.$
The proof is complete.
\endpf\\

Hence we have
\beq\notag
u=(x,y,-2z),\quad 
A=\left[
\begin{array}{rrr}
1&0&0\\
0&1&0\\
0&0&-2
\end{array}
\right],\qquad \Lambda={\bf 0}_{3\times 3}.
\eeq
\beq\notag
Ad=\big(\cos\theta\sin\va,\sin\theta\sin\va,-2\cos\va\big)^T
\eeq
\beq\notag
(d^TAd)d=(\sin^2\va-2\cos^2\va)\big(\cos\theta\sin\va,\sin\theta\sin\va,\cos\va\big)^T.
\eeq
Then 
\beq\notag
(d^TAd)d-Ad=-3\cos\va\sin\va\big(\cos\theta\cos\va,\sin\theta\cos\va,-\sin\va\big)^T,
\eeq

\bex
\Delta d+|\nabla d|^2d=\big(\va_{rr}+\frac{1}{r}\va_r-\frac{\sin(2\va)}{2r^2}\big)
\big(\cos\theta\cos\va,\sin\theta\cos\va,-\sin\va\big)^T,
\eex
\bex
d_t+u\cdot\nabla d=\left(\va_t+r\va_r\right)
(\cos\theta\cos\va,\sin\theta\cos\va,-\sin\va)^T.
\eex
Combining all above, the equation of $d$ in system \eqref{EL} becomes
\beq\label{ELaxeqn1}
\lambda_1(\va_t+r\va_r)=\va_{rr}+\frac{1}{r}\va_r-\frac{\sin(2\va)}{2r^2}-3\lambda_2\sin\va\cos\va
\eeq
If we assume $\lambda_2=0$, this is exactly the equation in \cite{HLLW16} for the simplified liquid crystal flows. In this paper, we would investigate the more general case when $\lambda_2\ne 0$.

The proof of the existence of a local smooth solution to equation \eqref{ELaxeqn1} proceeds similarly as in \cite{HLLW16} with minor modification to take care of the extra (smooth) external force term $3\lambda_2\sin\va\cos\va$. Therefore, we omit the proof here.

\subsection{Apriori estimates}
\label{secapri}

In general, we can not expect the energy inequality like \eqref{enginq} for the equation \eqref{ELaxeqn1} with the initial and boundary conditions \eqref{sbinitial} and \eqref{sboundary}. However, the following apriori estimate is still valid for \eqref{ELaxeqn1}.

\blm\label{axienglm}{\rm For any smooth solution $\va(r,t):[0,1]\times [0,T)\rightarrow \mathbb R$ to the equation \eqref{ELaxeqn1} with initial and boundary conditions \eqref{sbinitial} and \eqref{sboundary}, we have the following estimate
\beq
\int_0^1\left(|\va_r|^2+\frac{\sin^2\va}{r^2}\right)r\,dr+\int_0^T\int_0^1|\va_t|^2r\,dr\leq C(\va_0, T).
\eeq
}
\elm

\pf Multiplying the equation \eqref{ELaxeqn} by $\va_tr$ and integrating with respect to $r$ over $[0,1]$, we obtain 
\beq\notag
\begin{split}
&\frac{d}{dt}\frac12\int_0^1\left(|\va_r|^2+\frac{\sin^2\va}{r^2}\right)r\,dr+\lambda_1\int_0^1|\va_t|^2r\,dr\\
=&-\lambda_1\int_0^1\va_r\va_tr^2\,dr-3\lambda_2\int_0^1\sin\va\cos\va\va_tr\,dr
\\
\leq &\frac{\lambda_1}{2}\int_0^1|\va_t|^2r\,dr+C\int_0^1\left(|\va_r|^2+\frac{\sin^2\va}{r^2}\right)r\,dr.
\end{split}
\eeq
Standard Gronwall arguments should give us the desired energy estimates.

\endpf

We also need the following lemma.

\blm[{\rm comparison principle}]\label{comp}
{\it Suppose the functions $\va$, $f$ and $g$ are a smooth solution, subsolution and supersolution to \eqref{ELaxeqn1} on $[0,1]\times[0,T)$ respectively, and
$$f(r,t)\le \va(r,t)\le g(r,t) \ {\rm{on}}\ ([0,1]\times\{0\})\cup (\{0,1\}\times (0,T)).$$
Then we have
\beq
f(r,t)\le \va(r,t) \le g(r,t), \ \forall \ (r,t)\in [0,1]\times[0,T).
\eeq
}
\elm

\pf  Set $\bar f=f-\va$. Then $\bar f$ satisfies
\beq\notag%\label{lmpf2}
\displaystyle
\bar f_t+r\bar f_r\ge ~\bar f_{rr}+\frac{\bar f_r}{r}+p_2(r,t)\bar f,
\eeq
where
$$
\displaystyle p_2(r,t):=-\frac{\sin(2f(r,t))-\sin(2\va(r,t))}{2 r^2(f(r,t)-\va(r,t))}-3\lambda_2\frac{\sin(2f(r,t))-\sin(2\va(r,t))}{2(f(r,t)-\va(r,t))},
$$
and
$$
\bar f(r,t)\le 0 \ {\rm{on}}\ ([0,1]\times\{0\})\cup (\{0,1\}\times (0,T)).
$$
For any $t\in(0,T)$, there exists a small $r_2\in(0,1)$ such that $p_2(r,\tau)<0$ on $(r,\tau)\in(0,r_2)\times(0,t)$. Combining it with the fact that $p_2(r,\tau)$ is bounded on $(r_2, 1)\times(0,t)$, we conclude that $p_2(r,\tau)$ is bounded from above on $(r,\tau)\in(0,1)\times(0,t)$. By the standard maximum principle (see \cite{friedman} or \cite{pw}), 
we conclude that $\bar f\le 0$ or $\va\ge f$ on $[0,1]\times(0,T)$. 
Similarly, one can prove $\va\le g$.

\endpf
 
A natural corollary of the comparison lemma is the following maximum principle. 

\blm[{\rm maximum principle}]\label{maxp}
{\it Suppose that  $\va(r,t):[0,1]\times [0,T)\rightarrow \mathbb R$ is a smooth solution to equation \eqref{ELaxeqn1}.
If 
$$
\va_0(0)=0,  \  0\leq \va_0(r)\leq \pi, \ \forall \ r\in [0,1],
$$
then it holds
\beq
0< \va(r,t)<\pi, \ \forall\ r\in (0,1)\ {\rm{and}}\ t\in (0,T).
\eeq
}
\elm

\section{Existence and singularity of axisymmetric solutions}
 \setcounter{equation}{0}
\setcounter{theorem}{0}

In this section, we mainly concentrate on the equation
\beq\label{ELaxeqn}
\lambda_1(\va_t+r\va_r)=\va_{rr}+\frac{1}{r}\va_r-\frac{\sin(2\va)}{2r^2}-3\lambda_2\sin\va\cos\va
\eeq
 and provide the proof of Theorem \ref{glbth1} and Theorem \ref{blowthm2}.

\bigskip
\noindent{\bf Proof of Theorem \ref{glbth1}}.\quad
We argue by contradiction. 
Suppose $\va$ has the first singularity at finite time $0<T_1<\infty$. For any $r_0\in (0,1)$ and $t_1\in(0,T_1)$, by
the standard regularity theory of parabolic equations, we can prove
\beq\label{pf4}
\big\|\va\big\|_{C^k([r_0,1]\times(t_1,T_1))}\leq Cr_0^{-k}, \ \forall k\ge 1.
\eeq
This implies that the possible singularity of the  solution $\va$ can only happen at $r=0$. Combining this fact with the blowup criterion, Theorem 1.4 in \cite{huanglinwang14}, it holds for any $0<R<1$
\beq\label{blowupcr1}
\liminf\limits_{t\rightarrow T_1^-}\int_{0}^R|\va_r(r, t)|^2\,rdr\geq 4.
\eeq
We can choose $r_m\rightarrow 0^+$ and $t_1<t_m\rightarrow T_1^-$ as $m\rightarrow\infty$ and consider the blowup sequence 
\beq\label{rescaling}
\va^m(r,t)=\va\left(r_mr, t_m+r_m^2t\right).
\eeq
The rescaling functions $\va^m(r,t)$ are smooth solutions to the following equation
\beq\label{ELaxeqnre}
\lambda_1\left(\va^m_t+r_m^2\va^m_r\right)=\va^m_{rr}+\frac{1}{r}\va^m_r-\frac{\sin(2\va^m)}{2r^2}-3\lambda_2r_m^2\sin\va^m\cos\va^m,
\eeq
on $[0,1/r_m]\times[(t_1- t_m)/r_m^2,0]$. By the blowup criterion \eqref{blowupcr1}, we conclude that $\va^m$ is not constant.  
By the energy estimates in Lemma \ref{axienglm}, we have
\beq\label{reengest1}
\int_0^{1/r_m}|\va^m_r|^2\,rdr=\int_0^{1}|\va_r(r, t_m+r_m^2t)|^2\,rdr<\infty,
\eeq
and for any fixed $T>0$
\beq\label{reengest2}
\int_{-T}^0\int_0^{1/r_m}|\va^m_t|^2\,rdrdt=\int_{t_m-Tr_m^2}^{t_m}\int_0^{1}|\va_r(r,t)|^2\,rdrdt\rightarrow 0
\eeq
as $m\rightarrow \infty$. It is not hard to see that 
\beq\notag
r_m^2|\va^m_r|+3r_m^2|\lambda_2||\sin\va^m\cos\va^m|\rightarrow 0,\quad \mbox{in } L^2_{loc}([0,\infty)\times(-\infty,0))
\eeq
as $m\rightarrow \infty$. Therefore, there is a nonconstant smooth function $\omega(r):[0,\infty)\rightarrow \R$ with finite energy such that 
\beq\notag%\label{pf5}
\va^m(r,t)\rightarrow \omega\quad \mbox{in } C^2_{\mbox{loc}}([0,\infty)\times(-\infty,0]),
\eeq
and $\omega(r)$ is a solution to the axisymmetric  harmonic maps in dimension two
\beq\notag
\omega_{rr}+\frac{1}{r}\omega_r-\frac{\sin(2\omega)}{2r^2}=0.
\eeq

Combining the fact $0\leq\va_0< \pi$ with the lemma \ref{maxp}, we obtain
$$
0< \va(r,t)<\pi, \ r\in [0, 1], \ t_m\le t<T,
$$
which implies 
$$
0\leq\omega(r)\leq \pi\quad \mbox{for any }r\in[0,\infty).
$$
Hence, we conclude that 
$$
d_\infty=(\sin\omega\cos\theta, \sin\omega\sin\theta, \cos\omega):\R^2\rightarrow\mathbb S^2
$$  is a nontrivial harmonic map with finite energy. 

To secure a contradiction,
we need to construct suitable supersolutions and subsolutions to \eqref{ELaxeqn}. Denote
\beq\notag%\label{pf1}
\overline\va(r, t, c)=2\arctan \left(\frac{re^{bt}}{c}\right)
\quad\mbox{and}\quad
\underline\va(r,t,c)=2\arctan \left(-\frac{re^{bt}}{c}\right),
\eeq
for some positive constant $c$ and $b=3|\lambda_2|/\lambda_1$. It is easy to see that $\overline\va(r,c)$ and $\underline\va(r,c)$ are smooth functions. Direct calculation  implies
\beq\notag%\label{pf2}
\lambda_1(\overline\va_t+r\overline\va_r)-\overline\va_{rr}-\frac{\overline\va_r}{r}+\frac{\sin(2\overline\va)}{2r^2}+\frac{3\lambda_2}{2}\sin(2\overline\va)
=\frac{2rce^{bt}}{c^2+r^2e^{2bt}}\left(\lambda_1(1+b)+3\lambda_2\cos\overline\va\right)\geq 0,
\eeq
where we have used the facts $\lambda_1>0$ and $\left|\cos\overline\va\right|\leq 1$ in the last inequality.
Thus $\overline{\va}$ is a supersolution of \eqref{ELaxeqn}.
Denote
$$\eta_0=\pi-\max\limits_{0\leq r\leq 1}|\va_0(r)|.$$
By the assumption on $\varphi_0$, we have $\eta_0\ge 0$ and $|\va_0|\leq\pi-\eta_0$.
Since $\overline\va_r(0,0, c)=\frac{2}{c}$ and $\overline\va(0,t, c)=\va_0(0)=0$, we can find a sufficiently small $c>0$ 
such that
\beq\notag%\label{pf3}
\overline\va(r,0,c)\geq\va_0(r)
\eeq
 for any $r\in[0,1]$, with equality holds iff $r=0$. Similarly, we can prove that $\underline\va(r,t, c)$ is a subsolution to \eqref{ELaxeqn}
 and 
 \beq\notag%\label{pf3.1} 
 \varphi_0(r)\ge \underline\va(r,0,c) \ {\rm{in} } [0,1]
 \eeq
for a sufficiently small $c>0$. By lemma \ref{comp}, 
 {we conclude} that $\overline\va(r,t, c)\ge\va(r,t)\ge\underline\va(r,t, c)$ for $r\in [0,1]$ and $t>0$.
 
Now by the definition of $\overline\va(r,t, c), \underline\va(r,t, c)$, it is not hard to see that 
$$0< \va(r,t)<\pi-c_0, \ r\in [0, 1], \ t_m\le t<T,
$$
for some fixed constant $c_0>0$. Thus 
we have
$$
0\leq \omega(r)<\pi \quad \mbox{for any }r\in [0,\infty),
$$
which implies that  
$$
d_\infty=(\sin\omega\cos\theta, \sin\omega\sin\theta, \cos\omega):\R^2\rightarrow\mathbb S^2
$$ is a nontrivial harmonic map
with finite energy and degree zero, which is impossible.
Therefore, there is no finite time singularity for this case.
\endpf

\bigskip
Now, we will construct the suitable subsolution and initial data to show the existence of finite time singularity of \eqref{ELaxeqn}, by modifying the arguments in \cite{HLLW16} and \cite{chang-ding-ye} . 
\medskip

\noindent{\bf Proof of Theorem \ref{blowthm2}}.\quad We first consider the following potential subsolution to the equation \eqref{ELaxeqn}
\beq\notag
\eta(r,t)=2\arctan\left(\frac{r}{\beta(t)}\right),
\eeq 
where $\beta(t)$ satisfies the following ODE
\beq\label{eqbeta}
\frac{d\beta}{dt}=-\beta^{\frac23},\quad \beta(0)=\beta_0,
\eeq
where $\beta_0$ is a positive constant that will be determined later. It is not hard to see that 
\beq\label{beta}
\beta^{\frac13}(t)=\frac13(3\beta_0^{\frac13}-t)
\eeq
Setting $T_0=3\beta_0^{\frac13}$, we see that $\beta\rightarrow 0$ as $t\rightarrow T_0^-$. 
Direct computation implies 
\beq\label{eq4.8}
\begin{split}
&\lambda_1(\eta_t+r\eta_r)-\eta_{rr}-\frac{1}{r}\eta_r+\frac{\sin(2\eta)}{2r^2}+3\lambda_2\sin \eta\cos \eta\\
=&\frac{2r}{\beta^2+r^2}\left(\lambda_1\beta'+\lambda_1\beta+3\lambda_2\beta\cos\eta\right)\\
\leq &\frac{2r}{\beta^2+r^2}\left(-\lambda_1\beta^{\frac23}+(\lambda_1+3|\lambda_2|)\beta\right)\\
=&\frac{2r\beta^{\frac23}}{\beta^2+r^2}\left(-\lambda_1+(\lambda_1+3|\lambda_2|)\beta^{\frac13}\right).
\end{split}
\eeq
If we choose $\beta_0^{\frac13}<\frac{\lambda_1}{\lambda_1+3|\lambda_2|}$, then for any $(r,t)\in(0,1)\times(0,3\beta_0^{\frac13})$ it holds
\beq\label{subsol1}
\lambda_1(\eta_t+r\eta_r)-\eta_{rr}-\frac{1}{r}\eta_r+\frac{\sin(2\eta)}{2r^2}+3\lambda_2\sin \eta\cos \eta\leq 0.
\eeq
Since $\eta(1,0)< \pi$, we can find initial data $\va_0(r)$ such that $\va_0(0)=0$, $\va_0(1)>\pi$ and 
\beq\label{subsol2}
\eta(r,0)\leq \va_0(r)
\eeq
for any $r\in[0,1]$. By the definition of $\eta$, it is easy to see that for any $t\in(0,3\beta_0^{\frac13})$
\beq\label{subsol3}
\eta(0,t)=\va_0(0)=0, \quad \eta(1,t)\leq\pi< \va_0(1).
\eeq
Combining the facts \eqref{subsol1}-\eqref{subsol3} with Lemma \ref{comp}, we conclude that 
$
\eta(r,t)\leq \va(r,t)
$ 
for any $(r,t)\in[0,1]\times[0,3\beta_0^{\frac13})$. Direct computation implies 
\beq\notag
\eta_r(0,t)=\frac{2}{\beta(t)}\rightarrow \infty,\quad \mbox{as } t\rightarrow T_0^-.
\eeq
This combined with the fact $\eta(0,t)=\va(0,t)=0$ implies the finite time blowup of $\va(r,t)$, which completes the proof of Theorem \ref{blowthm2}.  
\endpf

\section{Second example of singularity}
\label{sec2ex}
 \setcounter{equation}{0}
\setcounter{theorem}{0}

We first need the following result on the existence of local smooth solution to (\ref{EL}) with the initial and boundary conditions \eqref{initial} and \eqref{bdry} under the assumptions \eqref{parodi}, \eqref{nec1} and \eqref{mus}.

\begin{theorem}\label{local-solution} For $(u_0, d_0)\in C^\infty_{0,{\rm{div}}}(B_1^3, \mathbb R^3)
 \times C_{\bf e}^\infty(B_1^3, \mathbb S^2)$, there exist $T_0 = T_0(u_0,d_0)> 0$
 and a unique smooth solution $(u,d) \in C^\infty\big(\overline{B_1^3} \times [0,T_0), \mathbb R^3 \times
 \mathbb S^2\big)$ to  the system  (\ref{EL}) 
 along with the initial-boundary conditions \eqref{initial} and \eqref{bdry}. 
 %under the assumptions \eqref{parodi}, \eqref{nec1} and \eqref{mus}. 
 Moreover, the energy dissipation inequality 
 \beq\label{energy_ineq}
\frac{d}{dt}\frac{1}{2}\int_{B_1^3}\left(|u|^2+|\nabla d|^2\right)\,
\leq-\int_{B_1^3}\left(\alpha_0|\nabla u|^2+\frac{1}{\lambda_1}\left|\Delta d+|\nabla d|^2d\right|^2\right)\,,
\eeq
holds for $0\le t<T_0$ and some $\alpha_0>0$.
\end{theorem} 

\pf The existence of local smooth solution has been proved in \cite{hieberpruss17}, and the energy inequality has been given by Lemma \ref{lemma0}.

\endpf

Now we would like to proceed to the proof of Theorem \ref{blowthm1}. Due to the energy estimate \eqref{enginq}, the proof is quite similar to that in \cite{HLLW16} for the simplified system. For the sake of completeness, we include a brief sketch of the main steps.

\bigskip 
\noindent{\bf Proof of Theorem \ref{blowthm1}}.  Assume $T_0 > 0$ is the end of the maximal time interval for the short time smooth solution $(u,d)$
from lemma \ref{local-solution}. Our goal is to show that 
If $\epsilon_0 >0$ is sufficiently small, then $T_0 <1$.

\smallskip
We argue by contradiction. Suppose for any $\epsilon>0$
we can find $(u_0, d_0) \in C^\infty_{0,{\rm{div}}} (B_1^3, \mathbb R^3) 
\times C^\infty_{\bf e}(B_1^3, \mathbb S^2)$ satisfying the assumption of Theorem \ref{blowthm1}, 
%such that $d_0$ is not homotopic to ${\bf e}$ relative to $\partial B_1^3$, and $E(u_0 , d_0) \le \epsilon^2$.
and a smooth solution $(u, d) \in C^\infty\big(\overline{B_1^3} \times  [0, 1], \mathbb R^3 \times\mathbb S^2\big)$
 to the Ericksen-Leslie system (\ref{EL})
with the initial-boundary conditions \eqref{initial} and \eqref{bdry} under the assumptions \eqref{parodi}, \eqref{nec1} and \eqref{mus}. Integrating (\ref{energy_ineq}) over $t$ yields that  
 $(u, d)$ satisfies the energy inequality
 \begin{equation}\label{energy_ineq1}
 E(u(t), d(t))+\int_0^t\int_{B_1^3} \big(\alpha_0|\nabla u|^2 + \frac{1}{\lambda_1}|\Delta d +|\nabla d|^2 d|^2\big)\le E(u_0, d_0) 
 \le\epsilon^2,
 \end{equation}
 for all $0\le t\le 1$.
Applying Fubini's theorem to (\ref{energy_ineq1}), we find that there exists $t_1 \in (\frac12, 1)$
such that
\begin{equation}\label{energy_ineq2}
 E(u(t_1), d(t_1))+\int_{B_1^3}\big(|\nabla u(t_1)|^2 + |\Delta d(t_1) + |\nabla d(t_1)|^2d(t_1)|^2\big)\le 8\epsilon^2.
 \end{equation}
From (\ref{energy_ineq2}) and $\epsilon$-apriori estimate Lemma 6.2 in \cite{HLLW16} on approximate harmonic maps, we conclude that there exists a constant $C > 0$ such that
 \begin{equation}\label{holder3}
 \big[d(t_1)\big]_{C^{\frac12}(B_1^3)} \le C\sqrt{\epsilon}.
 \end{equation} 
Thus $d(t_1)(B_1^3) \subset B^3_{C\sqrt{\epsilon}} ({\bf e}) \cap \mathbb S^2$
 and hence $d(t_1)$ is homotopic to $\bf e$ relative to $\partial B_1^3$, provided 
 $\epsilon >0$ is chosen to be sufficiently small. Since $d\in C^\infty\big(\overline{B_1^3}\times [0, t_1],\mathbb S^2\big)$
 and $d={\bf e}$ on $\partial B_1^3\times [0, t_1]$, we see that $d(t_1)$ is homotopic to $d_0$ relative to 
 $\partial B_1^3$ and hence $d_0$ is homotopic to $\bf e$ relative to 
 $\partial B_1^3$, which is a contradiction. This completes the proof of Theorem \ref{blowthm1}.
 \endpf

\section{Counter example on maximum principle}
\label{seccnex}
 \setcounter{equation}{0}
\setcounter{theorem}{0}

In this section, we will construct a counter example to show that the maximum principle is not necessarily valid for Ericksen-Leslie system. To this end, we consider the following Poiseuille flows for Ericksen-Leslie system \eqref{EL} in dimension one, i.e., for $(x,t)\in\mathbb R\times(0,\infty)$  
\begin{equation}\label{geneqn1D}
\begin{cases}
w_t+a=(g(\va)w_x+h(\va)\va_t)_x,\\
\lambda_1\va_t=\va_{xx}-h(\va)w_x,
\end{cases}
\end{equation}
where $g(\va)$ and $h(\va)$ are defined as following
\begin{align}\label{gh}\begin{split}
 g(\va):=&\mu_1\sin^2\va\cos^2\va+\frac{\mu_5-\mu_2}{2}\sin^2\va+\frac{\mu_3+\mu_6}{2}\cos^2\va+\frac{\mu_4}{2},\\
  h(\va):=&\mu_3\cos^2\va-\mu_2\sin^2\va=\frac{\lambda_1+\lambda_2\cos(2\va) }{2}.
 \end{split}
 \end{align}
 The system \eqref{geneqn1D} can be derived from the Ericksen-Leslie system \eqref{EL} by using the Poiseuille flows (see Appendix in \cite{CHL20} for details)
$$
u=(0,0,w(x,t)),\quad d=(\sin\varphi(x,t),0, \cos\varphi(x,t)).
$$
We consider a simplified case of \eqref{geneqn1D} by letting 
$$
a=0,\quad \quad \mu_1=\mu_5=\mu_6=0,\quad \mu_2=-1,\quad \mu_3=1, \quad \mu_4=3, \quad \lambda_1=\mu_3-\mu_2=2. 
$$
By the Onsager-Parodi relation (\ref{parodi}), one has
$$
\lambda_2=\mu_6-\mu_5=\mu_2+\mu_3=0,
$$
which implies 
$$
g(\va)= 2, \quad h(\va)=1.
$$
Then the system \eqref{geneqn1D} becomes
\begin{equation}\label{eqn1D}
\begin{cases}
w_t=2w_{xx}+\varphi_{tx},\\
2\varphi_t=\varphi_{xx}-w_x.
\end{cases}
\end{equation}
The energy inequality is as follows
$$
\frac{d}{dt}\frac12\int_{\R}|w|^2+|\varphi_x|^2\,dx+\int_{\R}|w_x|^2+|\varphi_t|^2+|w_x+\varphi_t|^2\,dx=0.
$$
Denote $v(x,t)=\int_{-\infty}^xw(y,t)\,dy$. Then it holds by the first equation of \eqref{eqn1D}
$$
v_t=2w_x+\varphi_t=2v_{xx}+\varphi_t.
$$
Therefore the system \eqref{eqn1D} can be written into a system of $v,\varphi$
\begin{equation}\label{eqn1Dv}
\begin{cases}
v_t=2v_{xx}+\varphi_{t},\\
2\varphi_t=\varphi_{xx}-v_{xx}.
\end{cases}
\end{equation}
Adding these two equations, we should have 
$$
(v+\varphi)_t=(v+\varphi)_{xx}.
$$ 
 Thus 
 $$
 (v+\varphi) (x,t)=H(x,t)*(v_0+\varphi_0)(x)
 =\frac{1}{\sqrt{4\pi t}}\int_{\R}\exp{\left(-\frac{|x-y|^2}{4t}\right)}(v_0+\va_0)(y)\,dy
 $$
 where $H(x,t)$ is the fundamental solution of heat equation in dimension one, and $(v_0, \varphi_0)$ is the initial data of the system \eqref{eqn1D}. Plugging in this back to the equation of $\varphi$, we should have 
 $$
 2\varphi_t=2\varphi_{xx}-\left(H(x,t)*(v_0+\varphi_0)(x)\right)_{xx}=2\varphi_{xx}-H(x,t)*(v_0+\varphi_0)_{xx}(x).
 $$
We may take the smooth initial data $v_0$ and $\va_0$ such that 
$$ 
v_0(x)=-x^2,\quad \va_0(x)=0.
$$
 Then
 $$H(x,t)*(v_0+\varphi_0)_{xx}=-2,$$ 
 and 
 $$
 \va(x,t)=t,
 $$ 
 which does not satisfy the maximum principle in the form of Lemma \ref{maxp}.

%\section{Blowup criteria in Serrin spaces}
% \setcounter{equation}{0}
%\setcounter{theorem}{0}
%
%
%
%In this section, we would extend the blowup criteria from simplified EL system to general EL system. The estimates for $d$ should be standard. The estimate for $u$ might be difficult, since $\sigma$ contains many second order terms and we need to make sure there are still parabolic. And also, there is a $d_t$ terms in equation of $u$, which makes this estimate even harder. 
%
%We may consider to extend the criteria of Huang-Wang CPDE to this case, which has been done by Wang-Zhang-Zhang for whole space. 

%\noindent{\bf Acknowledgments.}


\begin{thebibliography}{99}




\bibitem{chang-ding-ye} K. C.~Chang, W. Y.~Ding, R.~Ye,
Finite-time blow-up of the heat flow of harmonic maps from surfaces. J. Diff. Geom.,  36 (1992) 507-515.

\bibitem{CHL20}
G.~Chen, T.~Huang, and W.~S. Liu.
\newblock Poiseuille flow of nematic liquid crystals via the full
  Ericksen-Leslie model.
\newblock {\em Archive for Rational Mechanics and Analysis}, 236:839--891, 2020.



\bibitem{AnnaLiu19} F. de Anna and C. Liu, Non-isothermal General Ericksen–Leslie System: Derivation, Analysis and Thermodynamic Consistency. Arch. Ration. Mech. Anal., 231 (2019) 637-717.


\bibitem{DW} S. J. Ding, C. Y. Wang, Finite time singularity of the Landau-Lifshitz-Gilbert equation.
Int. Math. Res. Not. IMRN 2007, no. 4, Art. ID rnm012, 25 pp.


\bibitem{ericksen62} J. L. Ericksen,  Hydrostatic theory of liquid crystals. Arch. Ration. Mech. Anal. 9 (1962) 371-378.

\bibitem{friedman}
A. Friedman, Partial differential equations of parabolic type. Prentice-Hall, Englewood Cliffs, New Jersey, 1964.

\bibitem{hieberpruss17}
M. Hieber, J. Pr\"{u}ss, Dynamics of the Ericksen-Leslie equations with general Leslie stress I: the incompressible isotropic case. Math. Ann. 369(3-4) (2017) 977-996.

\bibitem{hongxin12} M.~C. Hong and Z.~P. Xin,
Global existence of solutions of the liquid crystal flow for the Oseen-Frank model in $\mathbb R^2$. Adv. Math. 231 (2012) 1364-1400.

\bibitem{huanglinwang14}
J. R. Huang, F. H. Lin and C. Y. Wang, 
Regularity and existence of global solutions to the Ericksen-Leslie system in $R^2$. Comm. Math. Phys., 331(2) (2014) 805-850. 

\bibitem{HLLW16}
T~Huang, F.~H. Lin, C.~Liu, and C.~Y. Wang,
Finite time singularity of the nematic liquid crystal flow in dimension three. Arch. Ration. Mech. Anal., 221(3) (2016) 1223--1254.

\bibitem{llwwz19}
C.~C. Lai, F.~H. Lin, C.~Y. Wang, J. C. Wei, and Y.~F. Zhou, Finite time blow-up for the nematic liquid crystal flow in dimension two.  Comm. Pure Appl. Math., (2021) https://doi.org/10.1002/cpa.21993.


\bibitem{leslie68} F. M. Leslie, Some thermal effects in cholesteric liquid crystals.
 Proc. Roy. Soc. A. 307 (1968) 359-372.



  \bibitem{Les79} F. M. Leslie, Theory of Flow Phenomena in Liquid Crystals. Advances in Liquid Crystals, Vol. 4, 1-81. Academic Press, New York, 1979. 
 
 
  \bibitem{LTX16} J. K. Li, E. Titi, and Z. P. Xin,
On the uniqueness of weak solutions to the Ericksen-Leslie liquid crystal model in $\mathbb R^2$. Math. Models Methods Appl. Sci. 26(4) (2016) 803-822.

 \bibitem{lin89} F. H. Lin, Nonlinear theory of defects in nematic liquid crystals: Phase transition and flow phenomena. Comm. Pure Appl. Math. 42 (1989) 789-814.

 \bibitem{linlinwang10} F. H. Lin, J. Y. Lin, and C. Y. Wang,
Liquid crystal flows in two dimensions. Arch. Ration. Mech. Anal. 197 (2010) 297-336.


\bibitem{linwang16} F. H. Lin and C. Y. Wang,
Global existence of weak solutions of the nematic liquid crystal
flow in dimension three.
Comm. Pure   Appl. Math. 69 (2016) 1532-1571.


\bibitem{linwang10} F. H. Lin and C. Y. Wang,
On the uniqueness of heat flow of harmonic maps and hydrodynamic
flow of nematic liquid crystals.
Chin. Ann. Math., Ser. B, 31 (2010) 921-938.

\bibitem{linwangbook}
F. H. Lin,  C. Y. Wang, The Analysis of Harmonic Maps and Their Heat Flows. The World Scientific, 2008.
 
  \bibitem{Parodi70} O. Parodi,
Stress tensor for a nematic liquid crystal. J. Phys. 31 (1970) 581-584.
 
\bibitem{pw}
M. H. Protter, H. F. Weinberge, Maximum principle in differential equations. Pretice-Hall, Englewood Cliffs, New Jersey, 1967.


\bibitem{wangwang14} M. Wang and W. D. Wang,
Global existence of weak solution for the 2-D Ericksen-Leslie system.
Calc. Var. Partial Differential Equations, 51(3-4) (2014) 915-962.

\bibitem{WZZ13} W. Wang, P. W. Zhang, and Z. F. Zhang,
Well-posedness of the Ericksen-Leslie system.
Arch. Ration. Mech. Anal. 210(3) (2013) 837-855.

\bibitem{WXL13}
H. Wu, X. Xu, C. Liu, On the general Ericksen-Leslie system: Parodi’s relation, well-posedness and stability. Arch. Ration. Mech. Anal. 208 (2013) 59-107.

 

\end{thebibliography}
\end{document}